\documentclass[a4paper]{jpconf}
\usepackage{graphicx}
\begin{document}
\title{Real-time plasma equilibrium reconstruction in a Tokamak.}

\author{J. Blum, C. Boulbe and B. Faugeras}

\address{Laboratoire J.A. Dieudonn\'e, UMR 6621, Universit\'e de Nice
  Sophia Antipolis, Parc Valrose, 06108 Nice Cedex 02, France}

\ead{jblum@unice.fr, boulbe@unice.fr, faugeras@unice.fr}

\begin{abstract}
The problem of equilibrium of a plasma in a Tokamak is a free boundary
problem described by the Grad-Shafranov equation in axisymmetric
configurations. The right hand side of this equation is a non linear
source, which represents the toroidal component of the plasma current
density. This paper deals with the real time identification of this
non linear source from experimental measurements.
 The proposed method is based on a fixed point algorithm, a finite
element resolution, a reduced basis method  and a least-square optimization formulation.  
\end{abstract}

\section{Introduction}
\begin{figure}
\begin{center}
\scalebox{0.36}{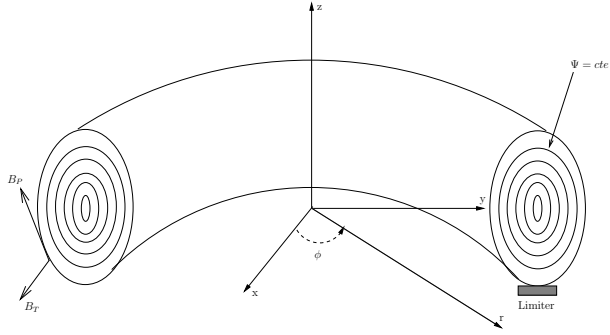}
\caption{Toroidal geometry.}\label{fig:geomtore}
\end{center}
\end{figure}

In Tokamaks, a magnetic field is used to confine a plasma in a
toroidal vacuum vessel. The magnetic
field is produced by external coils surrounding the vacuum vessel and
a current circulating in the plasma, forming a helicoidal resulting
magnetic field. At equilibrium, magnetic field lines lie on
isosurfaces forming a family of nested tori and called magnetic
surfaces. These magnetic surfaces enable to define the magnetic axis
and the plasma boundary. The innermost magnetic
surface which degenerates into a closed curve  is called magnetic axis. The plasma boundary
corresponds to the surface in contact with a limiter or being a
magnetic separatrix (hyperbolic line with an X-point).\\
 Denote by ${\bf j}$ the current density, ${\bf B}$ the magnetic field and $p$ the
kinetic pressure. 
The equilibrium of the plasma in presence of a magnetic field is
described by
\begin{eqnarray}
{\bf j}&=&\nabla\times\frac{{\bf B}}{\mu},\label{current}\\
{\bf j}\times{\bf B}&=&\nabla p,\label{equilibre}\\
\nabla\cdot{\bf B}&=&0,\label{div}
\end{eqnarray}
where $\mu$ is the magnetic permeability. Equation (\ref{current}) is
Ampere's theorem and equation (\ref{div}) represents the conservation of magnetic
induction.
Equation (\ref{equilibre}) means that at equilibrium the
Lorentz force ${\bf j}\times{\bf B}$ balances the force $\nabla p$
due to kinetic pressure. From this equation it is clear that
\begin{equation}
{\bf B}\cdot\nabla p=0 \mbox{ and }
{\bf j}\cdot\nabla p=0.
\end{equation}
Thus, $p$ is constant along magnetic field lines and current lines
which consequently lie on magnetic surfaces.\\
Consider the cylindrical coordinate system $(r,z,\phi)$. Under
axisymmetrical  hypothesis, the  magnetic field is supposed to be
independent of the toroidal angle $\phi$ and can be decomposed in the
form 
\begin{equation}\label{decomp}
\left\{
\begin{array}{lrc}
{\bf B}={\bf B}_P+{\bf B}_T\\
\displaystyle{{\bf B}_P}=\displaystyle{\frac{1}{r}[\nabla \psi\times {\bf e}_\phi]}\\
\displaystyle{{\bf B}_T}=\displaystyle{\frac{f}{r}{\bf e}_\phi}
\end{array}
\right.
\end{equation}
where ${\bf B}_P $ and ${\bf B}_T$ are respectively the poloidal and
the toroidal component of the magnetic field, ${\bf e}_\phi$ is a
unit vector in direction $\phi$,the term $\psi(r,z)$ is the poloidal
magnetic flux function (see Figure \ref{fig:geomtore}) and $f$ is the
poloidal current flux function. The magnetic surfaces are defined by rotation of
the flux lines $\psi=constant$ around the vertical axis $(Oz)$. 
Using (\ref{decomp}), equations (\ref{current})-(\ref{div}) lead to
the Grad-Shafranov equation (see \cite{Grad:1958,Shafranov:1958,Mercier:1974})
\begin{equation}
\label{eqn:gradshaf}
-\Delta^* \psi = r p'(\psi) + \frac{1}{\mu_0 r}(ff')(\psi)
\end{equation}
where
\begin{equation}
\label{eqn:delta*}
\Delta^* .=  \frac{\partial }{\partial r}(\frac{1}{\mu_0 r}  \frac{ \partial . }{\partial r}) 
+ \frac{\partial}{\partial z}(\frac{1}{\mu_0 r} \frac{ \partial .}{\partial z}).
\end{equation}
and $\mu_0$ is the magnetic permeability of the vacuum. 
Thus, under axisymmetric hypothesis, the three dimensional equilibrium
(\ref{current})-(\ref{div}) reduces to solve a two dimensional
non-linear problem. Note that the right hand side of
(\ref{eqn:gradshaf}) represents the toroidal component $j_{\phi}$ of the plasma
current density which is determinated by $p'$, $f$ and $f'$. \\
In this paper, we are interested in the numerical reconstruction of the
plasma current density and of the equilibrium (\ref{eqn:gradshaf}) (see \cite{Lao:1990,Blum:1990,Blum:1997}). This
reconstruction has to be achieved in real time  from experimental
measurements. The main difficulty consists in
identifying the functions $p'$ and $ff'$ in the non linear source term of
(\ref{eqn:gradshaf}). An iterative strategy involving a finite element
method to solve the direct problem (\ref{eqn:gradshaf}) and a least
square optimisation procedure to identify the non linearity using
reduced basis is proposed. A description of the experimental
measurements available in Tokamaks is given in Section 2. Section 3 is
devoted to the statement of the mathematical problem and the numerical
algorithm proposed. Numerical results obtained with the software Equinox are presented in Section 4.

\section{Experimental measurements}
\begin{figure}[!h]
\centering
\includegraphics[width=8cm,angle=-90]{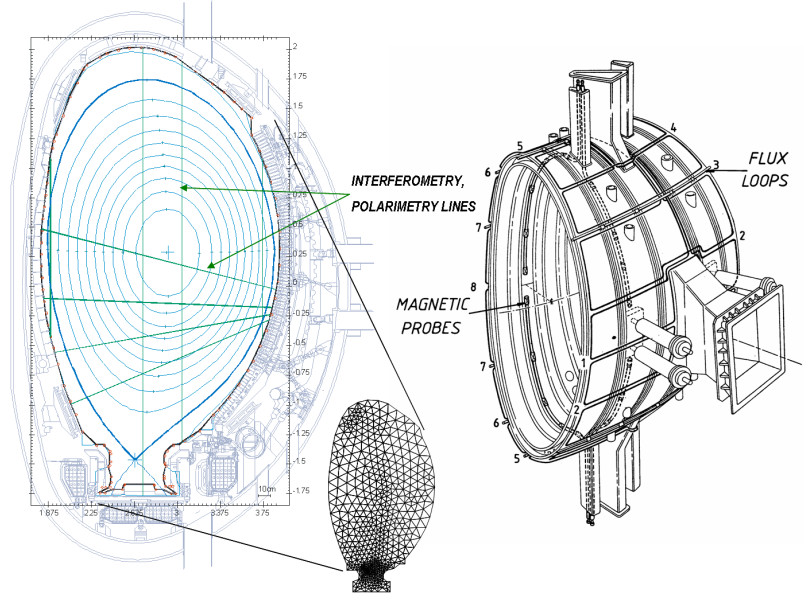} 
\caption{\label{fig:mes}  Left: the straight green lines represent the chords used for polarimetry and interferometry measurements. 
Right: part of the vacuum vessel including the magnetic measurements. 
At the bottom middle an example of finite element mesh used for numerical simulations.}
\end{figure}

Although $p$ and $f$ cannot be directly measured in a Tokamak, several
measurements are available (see Figure \ref{fig:mes}):
\begin{itemize}
\item magnetic measurements: the flux loops provide Dirichlet condition $\psi=h$
  and magnetic probes Neumann condition $\displaystyle{\frac{1}{r}\frac{\partial\psi}{\partial n}=g}$ 
  on the  walls of the vacuum vessel which represent the boundary of
  the computational domain. These measurements are given at some
  discrete points $N_i$ and the  Dirichlet full boundary data are obtained by interpolation. 
\item polarimetric measurements: this measurement gives the value of 
the integral along a family of chords $C_i$
$$
\int_{C_i} \frac{n_e}{r}\frac{\partial\psi}{\partial n}dl=\alpha_i.
$$
$n_e(\psi)$ is the electronic density which is approximately
constant on each flux line, $\displaystyle{\frac{\partial\psi}{\partial n}}$ is the
normal derivative of $\psi$ along the chord $C_i$.
\item interferometric measurements: they give the values of the
  integrals 
$
\displaystyle{\int_{C_i} n_e \,dl=\beta_i}.
$
\item current measurements: they give the value of the total plasma
  current $I_p$ defined by 
$$
I_p=\int_{\Omega_p} j_\phi dx.
$$
\end{itemize} 
Other measurements are potentially available but are not used for the moment in
the software Equinox.

\section{Algorithm and numerical resolution}
Let $\Omega$ be the domain representing the vacuum vessel of the
Tokamak, and $\partial \Omega$ its boundary. The equilibrium of a plasma in a Tokamak is a free boundary
problem. The plasma boundary is determinated either by its contact with a
limiter $D$ or as being a magnetic separatrix with an X-point (hyperbolic point). The region $\Omega_p\subset \Omega$ containing the plasma
is defined by
$$
\Omega_p=\{{\bf x}\in\Omega,\;\psi({\bf x})\geq \psi_b\}
$$
where $\psi_b=\max_D \psi$ in the limiter configuration or
$\psi_b=\psi(X)$ when an X-point exists.\\
In the vacuum region, the right hand side of (\ref{eqn:gradshaf}) vanishes
and the equilibrium reads
$$
\Delta^*\psi=0 \mbox{ in }\Omega\setminus\Omega_p
$$
Consider the following notations:
$\bar{\psi}= \displaystyle \frac{\psi - \displaystyle \max_\Omega \psi}{\psi_b - \displaystyle \max_\Omega \psi} \in [0,1]$ 
in $\Omega_p$ , 
$A(\bar{\psi})= \displaystyle \frac{R_0}{\lambda}p'({\psi})$ and 
$B(\bar{\psi})=\displaystyle \frac{1}{\lambda\mu_0 R_0}(ff')({\psi})$. 
The functions $A$ and $B$ are flux functions defined on the
fixed interval $[0,1]$.
Giving
Dirichlet boundary conditions, the final 
equilibrium problem is
\begin{equation}
\label{eqn:final}
\left \lbrace
\begin{array}{rcl}
-\Delta^* \psi &= &\lambda[\displaystyle\frac{r}{R_0} A(\bar{\psi}) +  \displaystyle 
\frac{R_0}{r} B(\bar{\psi})] \chi_{\Omega_p}\quad \mathrm{in}\ \Omega \\[10pt]
\psi&=& h\quad  \mathrm{on}\  \partial \Omega 
\end{array} 
\right.
\end{equation}
where  $\chi_{\Omega_p}$ 
is the characteristic function of $\Omega_p$ and $R_0$ is the major
radius of the Tokamak.
The parameter $\lambda$
is a normalization factor that satisfies
\begin{equation}\label{lambda}
I_p=\lambda \int_{\Omega_p}[\displaystyle\frac{r}{R_0} A(\bar{\psi}) +  \displaystyle 
\frac{R_0}{r} B(\bar{\psi})]  d\Omega.
\end{equation}

\subsection{Iterative algorithm}
The aim of the method is to reconstruct the plasma current
density and the equilibrium solution in real time. At each time step
determined by the availability of new measurements, the method
consists in constructing a sequence $(\psi^n, \Omega_p^n, A^n, B^n)$ converging to the solution vector $(\psi, \Omega_p, A, B)$. 
The sequence is obtained by the following algorithm:
\begin{itemize}
\item Starting guess: $\psi^0$, $\Omega_p^0$, $A^0$ and $B^0$ known from the previous time step solution. 
Compute $\lambda$ satisfying (\ref{lambda}).
    \item Step 1 - Optimisation step: computation of
      $A^{n+1}(\bar{\psi}^n)$, $B^{n+1}(\bar{\psi}^n)$ using a
      least square procedure.
   \item Step 2 - Direct problem step: computation of $\psi^{n+1}$ and
     $\Omega_p^{n+1}$ solution of
\begin{equation}\label{GSA}
\left \lbrace
\begin{array}{rcl}
-\Delta^* \psi^{n+1} &= &\displaystyle{\lambda [\frac{r}{R_0} A^{n+1}(\bar{\psi}^{n})} +  \displaystyle 
\frac{R_0}{r} B^{n+1}(\bar{\psi}^{n})] \chi_{\Omega_p^n}\quad \mathrm{in}\ \Omega \\[10pt]
\psi^{n+1}&=& h\quad  \mathrm{on}\  \partial \Omega. 
\end{array} 
\right.
\end{equation}

\item Step 3: If the process has not
  converged, $n:=n+1$ and return to Step 1 else $(\psi, A, B)=(\psi^n, A^n, B^n)$. 
 The process is supposed to have converged when the
  relative residu
  $\displaystyle{\frac{||\psi^{n+1}-\psi^{n}||}{||\psi^{n}||}}$ is
  small enough.
 \end{itemize}
At each iteration of the algorithm, an
inverse problem corresponding to the optimization step and an 
approximated direct Grad-Shafranov problem corresponding to (\ref{GSA})
 have to be solved successively. In (\ref{GSA}), $\psi^n$ is
 supposed known, the right hand side does not depend on
 $\psi^{n+1}$ and thus problem (\ref{GSA}) is linear.\\
Step 1 consists in a least-square minimization formulation 
\begin{equation}
\left \lbrace
\begin{array}{l}
\mathrm{Find}\ A^*,\ B^*,\ n_e^*\ \mathrm{such}\ \mathrm{that}: \\[10pt]
J(A^*,B^*,n_{e}^*)=\inf J(A,B,n_{e}). 
\end{array}
\right.
\end{equation}
 When polarimetric and interferometric measurements are used, the electronic density has also to be identified even if $n_e$, depending on $\bar\psi$, does not appear in (\ref{eqn:final}).
The cost function $J$ is defined by
\begin{equation}
\label{eqn:costfunction}
J(A,B,n_e)=J_0 + K_1 J_1 + K_2 J_2 + J_\epsilon
\end{equation}
where
$$
 J_0=  \displaystyle \sum_i 
 ( \displaystyle \frac{1}{r} \frac{\partial \psi}{\partial n}(N_i) 
 - g_i)^2 ,\;\; 
 J_1=  \displaystyle \sum_i (\displaystyle \int_{C_i} 
 \frac{n_e}{r} \frac{\partial \psi}{\partial n}dl - \alpha_i)^2,
 \mbox{ and  }
J_2=\displaystyle \sum_i (\displaystyle \int_{C_i} n_edl - \beta_i)^2,\;\;
$$
and the coefficients $K_1$ and $K_2$ are weights giving more or less
importance to measurements used (see  \cite{Blum:1990}). As a consequence of the ill-posedness of the identification of $A$, $B$ and $n_e$, a Tikhonov regularization term $J_\epsilon$ is introduced (see \cite{Tikhonov:1977}) where
$$
J_\epsilon = 
\epsilon_1 \displaystyle \int_0^1 [A''(x)]^2 dx 
+
\epsilon_2 \displaystyle \int_0^1 [B''(x)]^2 dx 
+
\epsilon_3 \displaystyle \int_0^1 [n_e''(x)]^2 dx 
$$
and $\epsilon_1$, $\epsilon_2$ and $\epsilon_3$ are regularization parameters.\\
In the next section, the numerical methods used to solve the two involved
 problems are detailed.

\subsection{Numerical resolution}
\label{sec:numerical-resolution}
The resolution of the direct problem is based on a $P^1$ finite element
method (see \cite{Ciarlet:1980}). Consider the family of triangulation ${\tau}_h$ of
$\Omega$, and $V_h$ the finite dimensional subspace of $H^1(\Omega)$
defined by
$$
V_h=\{v_h\in H^1(\Omega), v_{h|T}\in P^1(T),\,\forall T\in {\tau}_h\}.
$$
Introduce $V_h^0=V_h\cap H^1_0(\Omega)$, the discrete variational
formulation of (\ref{GSA}) reads
\begin{equation}
\label{eqn:var}
\left \lbrace
\begin{array}{l}
\mathrm{Find}\ \psi_h\in V_h \mbox{ with }
\psi_h=h\mbox{ on }\partial\Omega \mbox{ such that }\\[10pt]
\displaystyle \forall v_h \in V_h^0, \int_\Omega \displaystyle \frac{1}{\mu_0 r}\nabla \psi_{h} \cdot \nabla v_h dx = 
\int_{\Omega_p}\lambda [\displaystyle \frac{r}{R_0} A(\bar{\psi^*})+
  \displaystyle \frac{R_0}{r}B(\bar{\psi^*})] v_h
dx \\[10pt]
\end{array}
\right .
\end{equation}
where $\psi^*$ represents the value of $\psi$ at the previous iteration.
The Dirichlet boundary conditions are imposed using the method
consisting in computing the stiffness matrix of the Neumann problem and
modifying it. The modifications consist in replacing the rows corresponding
to each boundary node setting $1$ on
the diagonal terms and $0$ elsewhere. The Dirichlet conditions appear
in the right hand side of the linear system. At each
iteration, only the right hand side of the system has to be modified. 
 The sparse matrix $K$ and its inverse $K^{-1}$ are computed at the beginning of the
algorithm and stored until the end of the simulation. \\
The linear system of (\ref{eqn:var}) can be written in the form
\begin{equation}\label{linear}
K.\Psi=y+H
\end{equation}
where $K$ is the modified stiffness matrix, $\Psi$ is the unknown vector, $y$ is the right hand side of the
problem and $H$ is the term corresponding to the Dirichlet
conditions. The vector $y$ depends on $A$ and $B$ determinated in
the optimization step. \\
The functions $A$, $B$ and $n_e$ are decomposed on a basis $(\Phi_i)_{i=1,...,m}$ of dimension $m$ 
$$A(x)=\sum_i^m a_i\Phi_i(x),\;\;B(x)=\sum_i^m b_i\Phi_i(x)\mbox{ and
} n_e(x)=\sum_i^m c_i\Phi_i(x).$$

The vector $y$ reads
\begin{equation}\label{linearRhs}
y=Y(\bar{\psi^*})u
\end{equation}
where $u=(a_1,...,a_m,b_1,...,b_m)\in\mathbf{R}^{2m}$. The term $Y$ is a matrix of size $n\times 2m$ where $n$ is the number of nodes. Consider $(v_i)$ a basis of $V_h$, each row $i$ of $Y$ is decomposed as
$$
Y_{ij}(\bar\Psi^*)=\left\{
\begin{array}{llc}
\displaystyle{\int_{\Omega_p} \lambda \displaystyle\frac{r}{R_0}\Phi_j(\bar{\psi^*}) v_i d\Omega} &\mbox{ if }1\leq j\leq m\\
\displaystyle{\int_{\Omega_p} \lambda \displaystyle\frac{R_0}{r}\Phi_{j-m}(\bar{\psi^*}) v_i d\Omega} &\mbox{ if }m+1\leq j\leq 2m.
\end{array}
\right.
$$
During the optimisation step, $n_e$ is first estimated from
interferometric measurements and $A$ and $B$ are computed in a second time.
The function $n_e(\bar{\psi})$ is approximated using a least square
formulation for the minimun of $J_2$ with Tikhonov regularization, and solving the associated normal equation.\\
To approximate $A$ and $B$, suppose $n_e$ is known and consider the discrete approximated inverse problem
\begin{equation}
\label{eqn:optim}
\left \lbrace
\begin{array}{l}
\mathrm{Find}\ u\ \mathrm{minimizing}:\\
J(u)=\displaystyle\frac{1}{2}
\|  C({\psi^*},n_e)\Psi - k \|^2_D 
+\displaystyle \frac{\varepsilon}{2}u^{T} \Lambda u
\end{array}
\right.
\end{equation}
where $C(\psi^*,n_e)$ is the observation operator corresponding to
experimental measurements given in Section 2 and $k$ represents the experimental measurements. 
The first term in $J$ is the discrete version of $J_0+K_1 J_1$. 
The second one corresponds to the first two terms of the
Tikhonov regularization term $J_\varepsilon$. 
Denote by $l$ the number of measurements available and by $\sigma_i^2$ the variance 
of the error associated to the $i^{th}$ measurement, the norm $\|.\|_D$ is defined by
$$
\forall {\bf x}\in\mathbf{R}^l\;\; \|{\bf x}\|_D^2=(D {\bf x},{\bf x})=(D^{1/2}{\bf x},D^{1/2}{\bf x})
$$
where $D$ is the diagonal matrix $\displaystyle d_{ii}=\frac{1}{\sigma_i^2}$.\\
$C$ is a matrix of size $l\times n$ and can be viewed as a vector
composed of two blocks corresponding respectively to $J_0$ and $J_1$. \\
The matrix $\Lambda$ is of size $2m\times 2m$ and is block diagonal composed of two blocks $\Lambda_1$ and $\Lambda_2$ of size $m\times m$, with
$$
(\Lambda_1)_{ij}=(\Lambda_2)_{ij}=\int_0^1 \Phi"_i(x)\Phi"_j(x)dx
$$
where the $\Phi^"_i$ are the second derivatives of the basis functions $\Phi_i$. \\
Using (\ref{linear}) and (\ref{linearRhs}), the problem (\ref{eqn:optim}) becomes
$$
\begin{array}{lll}
J(u)&=&\displaystyle\frac{1}{2}\|  C({\psi^*},n_e)\Psi - k \|_D^2 + \displaystyle\frac{\varepsilon}{2}u^{T} \Lambda u \\[10pt]
&=&\displaystyle\frac{1}{2}\|  C({\psi^*},n_e) {K}^{-1} Y(\bar{\psi^*}) u +( C({\psi^*},n_e){K}^{-1} H - k ) \|^2_D +\displaystyle\frac{\varepsilon}{2} u^{T} \Lambda u \\[10pt]
&=&\displaystyle\frac{1}{2}\|  E u- F \|^2_D + \displaystyle\frac{\varepsilon}{2}u^{T} \Lambda u 
\end{array}
$$
where $E= C({\psi^*},n_e) {K}^{-1} Y(\bar{\psi^*})$ and $F= -C({\psi^*},n_e){K}^{-1} H + k$.
Setting $\tilde{E}=D^{1/2}E$, problem (\ref{eqn:optim}) reduces to solve the normal equation
\begin{equation}
\label{eqn:normal}
(\tilde{E}^T \tilde{E}+\varepsilon\Lambda)u=\tilde{E}^T F
\end{equation}

\section{Numerical results}
The method detailed in this paper has been implemented in the software
Equinox developed in collaboration with the Fusion Department at
Cadarache for Tore Supra (see \cite{Bosak:2001,Blum:2004,Bosak:2003}) 
and JET (Join European Torus). Equinox can be used on the one hand for precise studies in which the computing time is not a limiting factor 
and on the other hand in a real-time framework.

In table \ref{tab:residu}, the evolution of the relative
residu on $\psi$, $A$ and $B$ versus the number of iterations is given. 
It demonstrates numerically the convergence of the algorithm. 

\begin{table}
\caption{\label{tab:residu} Numerical convergence.}
\begin{center}
\begin{tabular}{llll}
\br
Iteration $n$ &$\displaystyle{\frac{\|\psi^{n+1}-\psi^n\|}{\|\psi^n\|}}$ &$\displaystyle{\frac{\|A^{n+1}-A^n\|}{\|A^n\|}}$&$\displaystyle{\frac{\|B^{n+1}-B^n\|}{\|B^n\|}}$\\
\mr
1&0.0001477 &0.0161846 &0.0212656\\
2&1.77698e-05 &0.0008774 &0.0021739\\
3&6.74457e-06 &0.0001095& 0.0003792\\
4&4.71623e-06 &1.1026e-05 & 7.0184e-05\\
5&3.45322e-06& 1.70946e-05& 2.07231e-05\\
6&2.53969e-06 &1.27861e-05& 1.25182e-05\\
7&1.86885e-06& 8.90356e-06 &9.60398e-06\\
8&1.37495e-06 &6.29852e-06& 7.40854e-06\\
9&1.01137e-06& 4.54006e-06& 5.59427e-06\\
10&7.43857e-07& 3.30811e-06 &4.1659e-06\\
\br
\end{tabular}
\end{center}
\end{table}

If a value $\epsilon=10^{-6}$ is used as stop condition, 
the algorithm almost always needs about $5$ to $30$ iterations to converge. 
For the finite element mesh of $412$ nodes and $762$ elements used at Tore Supra 
the average iteration cost (file access included) 
is of $0.038\ s$ on a processor Intel(R) Core(Tm)2 Duo T7100 1.80GHz.
The expensive operations are the updates of matrices $C$ and $Y$ and the computation of products $CK^{-1}$ and $CK^{-1}Y$ 
(see Section \ref{sec:numerical-resolution}). The resolution of the direct problem (Eq. \ref{linear}) is cheap since the 
matrix $K^{-1}$ is stored once for all and does not change from one iteration to the other. 
Morevover the resolution of the normal equation (Eq. \ref{eqn:normal}) is also cheap. Each function $A$ and $B$ is decomposed in a basis 
$\Phi_i$ (which can be chosen to be cubic B-splines, piecewise linear functions or wavelet scaling functions) 
of about $5$ to $10$ functions (depending on the user's choice). This makes 
the dimension of the matrix to be inverted at most $20 \times 20$.

For real-time applications the numerical reconstruction of an equilibrium must not 
take more than $0.050$ to $0.100\ s$.  
Therefore in this context it is not possible to let the algorithm fully converge.
A maximum number of iterations has to be set 
(typically $1$ or $2$ depending on the computer and the size of the mesh). 
In practice this is not a real problem since our approach is quasi-static. Two successive equilibriums 
during a pulse are very close and no important differences have been observed between the results of the fully converged algorithm 
and its real-time version.

Concerning Tikhonov regularization with a real-time constraint it is absolutely not possible 
to compute the regularization parameter for each equilibrium. 
Therefore it is kept constant during a pulse. The L-curve method \cite{Hansen:1998} provided a value of approximately $5\times10^{-5}$.

\begin{figure}[!h]
\centering
\includegraphics[width=14cm,height=10cm,angle=0]{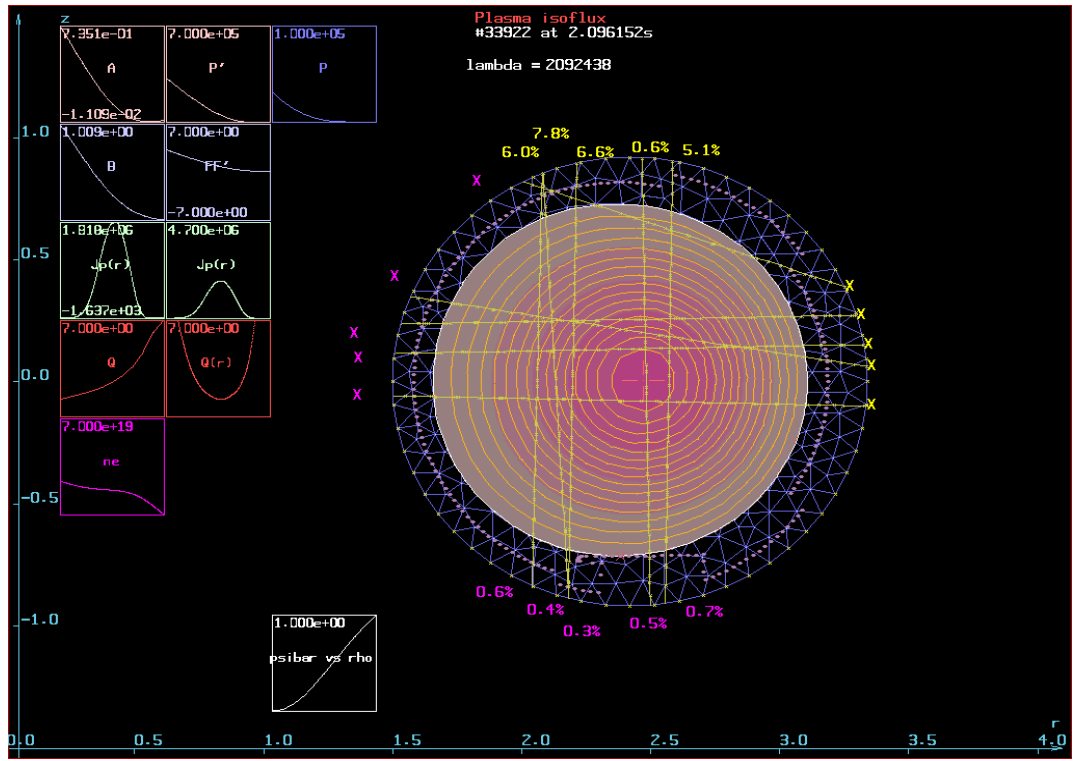} 
\caption{\label{fig:magpolar} Example of Equinox output. See text for details.}
\end{figure}

An example of the outputs of Equinox is presented in figure \ref{fig:magpolar}. 
It is a Tore Supra pulse in which magnetic, interferometric and polarimetric measurements on 5 chords are used. 
One can observe the position of the plasma in the vacuum vessel. 
Isoflux lines are displayed from the magnetic axis to the boundary. 
For each active chord, the error between computed and measured interferometry is given in purple. 
These errors are less than $1\%$. The polarimetry errors are given in yellow. 
They are one order of magnitude larger than for interferometry. 
Different graphs are plotted on the left of the display. 
On the first row the identified function $A$, and corresponding functions $p'$ and $p$. On the second row the identified function $B$ and corresponding function $ff'$. 
The third row gives the toroidal current density $j_{\phi}$ in the equatorial plane  
and the fourth one shows the safety factor. 
Finally on the fifth row the identified $n_e$ function is plotted.

\section*{References}
\bibliographystyle{iopart-num}
\bibliography{biblio-fusion}

\end{document}